\documentclass{elsart}

 \newtheorem{Corollary}{Corollary}
 \newtheorem{Theorem}{Theorem}
 \newtheorem{Proposition}{Proposition}
 \newtheorem{Lemma}{Lemma}
  \newcommand{\epf}{\hfill$\Box$}
  \newcommand{\su}[1]{\mathrm{Supp}\,#1}

\usepackage{amsfonts}
 
\begin{document}

 \begin{frontmatter}



\title{Finite dimensional graded simple algebras}


\author[yb]{Y. A. Bahturin}
\address{Department of Mathematics and Statistics\\Memorial
University of Newfoundland\\ St. John's, NL, A1C5S7,
Canada and \\ Department of
Algebra, Faculty of Mathematics and Mechanics\\Moscow
State University\\Moscow, 119992, Russia}
\ead{yuri@math.mun.ca}\thanks[yb]{Work is partially supported by NSERC grant \#
227060-04 and URP grant, Memorial University of Newfoundland}

\author[sks]{S.K.Sehgal}
\address{Department of Mathematical Sciences\\University of Alberta\\ Edmonton, AB\quad T6G 2G1, Canada}
\ead{S.Sehgal@math.ualberta.ca}\thanks[sks]{Work is partially supported by NSERC grant \#
A-5300}

\author[mz]{M. V. Zaicev}
\address{Department of Algebra,
Faculty of Mathematics and Mechanics\\Moscow State
University\\Moscow, 119992, Russia}\ead{zaicev@mech.math.msu.su}
\thanks[mz]{Work is partially supported by
RFBR, grant 02-01-00219, and SSC-1910.2003.1}


\vskip 0.5in
\begin{abstract}
Let $R$ be a finite-dimensional algebra over an algebraically closed field
$F$ graded by an arbitrary group $G$. We prove that $R$ is a graded division
algebra if and only if it is isomorphic to a twisted group
algebra of some finite subgroup of $G$. If the characteristic of $F$ is zero
or ${\rm char}~F$ does not divide the order of any finite subgroup of $G$
then we prove that $R$ is graded simple if and only if it is a matrix
algebra over a finite-dimensional graded division algebra.
\end{abstract}
\end{frontmatter}

\section{Introduction}\label{s1}

Graded rings and graded algebras have been investigated intensively during the last few 
decades. One of the important directions of this investigation is the 
development of the structure theory of graded rings and algebras, in particular
the description of radicals, simple and semisimple objects. One of the first
steps on this way is the description of all possible gradings on simple 
algebras. For example, in \cite{K} all integer gradings on finite 
dimensional simple complex Lie algebras have been described. In \cite{Kac} all 
gradings by a finite cyclic group on complex simple finite-dimensional Lie
algebras have obtained. In \cite{E} all possible gradings on the
Cayley-Dickson algebra have been classified. In \cite{BaShest} and \cite{BShZ}, \cite {BZA}, and \cite{BSh}
authors describe all Abelian gradings on classical simple finite-dimensional 
simple Lie and special Jordan algebras.

During the last few years  a number of papers   have been published, dedicated group grading on simple
associative algebras. In \cite{Smirn} the author
has shown that finite ${\bf Z}$-gradings on a simple associative algebra,
not necessary finite-dimensional, can be obtained from the Peirce
decomposition. In \cite{DINR} all gradings of matrix algebra by a torsion
free group have been described. The latter result was generalized in 
\cite{ZS} to any simple Artinian ring. In \cite{BSegZ} all Abelian gradings 
on matrix algebra over an algebraically closed field of characteristic
zero have been precisely  described. Later in \cite{BZ} this result was generalized to non-abelian gradings on matrix algebras.

It is well-known that a graded simple algebra is not in general simple in the
ordinary sense. For example, the group algebra $F[G]$ of any non-trivial
group over a field is not simple but is
graded simple in the canonical $G$-grading. So, the description of all 
possible gradings on simple algebras is only the first step toward the 
classification of graded simple algebras.

One of the first results in this area became the classification of finite
dimensional associative simple superalgebras, i.e. associative algebras
with ${\bf Z}_2$-gradings, over an algebraically closed field \cite{W}.
Now this result is well - under\-stood and any such  superalgebra is
either a matrix algebra with a ${\bf Z}_2$-grading or tensor product 
of a matrix algebra with the group algebra of ${\bf Z}_2$. The progress in the
description of all abelian gradings on full matrix algebra made it possible to complete the description of all graded simple algebras. In
already quoted paper \cite{BSegZ} the authors classified finite-dimensional graded 
simple algebras over an algebraically closed field of characteristic zero 
provided that the grading group is abelian. In fact, for a wide class of
abelian groups the result looks like in the super-case: $G$-graded algebra 
$A$ is graded simple if and only if $A$ is the tensor product of a matrix
algebra with some $G$-grading with the group algebra of some finite 
subgroup of $G$. Some partial results about graded simple algebras
with non-abelian torsion free grading have also been obtained in \cite{BSegZ}.

In the present paper we study finite-dimensional graded simple algebras over 
an algebraically closed field $F$ under some natural restriction on the
characteristic of $F$. In Section \ref{s2} we recall all necessary definitions
and notions. In Section \ref{s3} we prove the unitarity of graded simple algebras.
Main result of Section \ref{s4} is that any graded division algebra is a twisted
group algebra of some finite group $G$. In Section \ref{s5} we prove the 
semiprimitivity of graded simple algebra under a weak assumption on
${\rm char}~F$. The main result of this paper, presented in Section \ref{s6},
Theorem \ref{T3} states that any graded simple algebra is a matrix algebra over
a graded division algebra under the same restriction on the characteristic
of base field as in the Section \ref{s5}.

\section{Preliminaries}\label{s2}

Here we recall main definitions and constructions. Let $R$ be an associative
algebra over a field $F$ and $G$ be a multiplicative group with the
identity element $e$. Then $R$ is said to be a $G$-{\it graded algebra} (or simply {\it
graded} if $G$ is fixed) if there is a vector space decomposition
$$
R=\bigoplus_{g\in G} R_g
$$
where the subspaces $R_g$ satisfy $R_gR_h\subseteq R_{gh}$ for all $g,h\in G$.
The subspaces $R_g$ are called the {\it homogeneous components} of the grading,
in particular, $R_e$ is called the {\it identity component}. Clearly,
$R_e$ is a subalgebra of $R$. Given $g\in G$, any element $a\in R_g$ is
called {\it homogeneous} and we write $\deg a=g$. By definition of the grading,
any $x\in R$ can be uniquely written as
$$
x=\sum_{g\in G} x_g \qquad{\rm with}~~x_g\in R_g.
$$
A subspace $V\subseteq R$ is called a {\it graded subspace} (or
{\it homogeneous subspace}) if
$$
V=\bigoplus_{g\in G} (V\cap R_g).
$$
Equivalently, $V$ is graded if and only if for any $x=\sum_{g\in G} x_g\in V$ we have
$x_g\in V$ for any $x\in V$and $g\in G$. A 
subalgebra (left-, right-, two-sided ideal) $A$ is said to be {\it graded}
if it is a graded subspace. Following this general approach we say that
$R=\oplus_{g\in G}R_g$ is a {\it graded simple algebra} if $R^2\ne 0$ and $R$
has no non-trivial proper graded ideals. Similarly, $R$ is a {\it graded
division algebra} if $R$ is unitary, that is, possesses a multiplicative unit $1$, and any non-zero homogeneous element
of $R$ is invertible.

Let $A=\oplus_{g\in G}A_g$ and $B=\oplus_{g\in G}B_g$ be two $G$-graded
algebras. An algebra homomorphism (isomorphism) $f:A\rightarrow B$ is called a
{\it homomorphism (isomorphism) of graded algebras} if $f$ preserves the 
graded structure, that is, $f(A_g)\subseteq B_g$ for all $g\in G$.

Given an arbitrary algebra $R$, the easiest way to define a $G$-grading on
$R$ is to set $R_e=R$ and $R_g=0$ for all $e\ne g\in G$. This is the so-called
{\it trivial} $G$-grading. Thus, in a
$G$-graded algebra $R=\oplus_{g\in G}R_g$ it is not necessary that all homogeneous
components $R_g$ are non-zero. The set 
$$
{\rm Supp}~R=\{g\in G\vert R_g\ne 0 \}
$$
is called the {\it support} of the grading.
In general the support of a $G$-grading is not a subgroup of $G$. For
example, if $R^2=0$ then ${\rm Supp}~R$ could be an arbitrary subset of the grading group.
As we will see later, ${\rm Supp}~R$ need not be a 
subgroup, even if the multiplication in $R$ is non-trivial.

As an illustration of the preceding definitions we present few examples of graded
algebras which are essential for our work.

\vskip .2in

\subsection{Elementary grading on full matrix algebras}\label{ss1}
 Let 
$R=M_n(F)$ be a $n\times n$-matrix algebra over $F$ and $G$ a group. 
Fix an arbitrary $n$-tuple ${\bf g}=(g_1,\ldots, g_n)\in G^n$ of elements 
of $G$. Then ${\bf g}$ defines a $G$-grading on $R$ as follows. Let
 $E_{ij}, 1\le i,j\le n$, be the set of all matrix units of $R$.  Then we set
$$
R_g=\mathrm{Span}\,\{E_{ij}\vert g_i^{-1}g_j=g \}.
$$
Direct computations show that $R_gR_h\subseteq R_{gh}$ for all $g,h\in G$
and hence the decomposition
$$
M_n(F)=R=\bigoplus_{g\in G} R_g
$$
is a $G$-grading. We call this grading an {\it elementary grading} defined
by the $n$-tuple ${\bf g}$.

Note that if ${\bf g'}=(ag_1,\ldots, ag_n)$ is another $n$-tuple with
an arbitrary $a\in G$ then it defines the same $G$-grading since
$(ag_i)^{-1}(ag_j)=g_i^{-1}g_j$. In particular, one can always assume $g_1$ being the identity element of $G$. Given any permutation $\sigma$ of integers
$\{1,\ldots,n\}$, we may consider ${\bf g}_\sigma=(g_{\sigma(1)},\ldots,
g_{\sigma(n)})$. Then ${\bf g}$ and ${\bf g}_\sigma$ define on $R=M_n(F)$
the structure of two $G$-graded algebras $R({\bf g})$ and 
$R({\bf g}_\sigma)$. But they are isomorphic as $G$-graded algebras. This
can be viewed, for instance, in the following way.

Let $V$ be a $n$-dimensional vector space over $F$. We say that $V$ is
a $G$-{\it graded space} if $V=\oplus_{g\in G} V_g$. A linear transformation
$f:V\rightarrow V$ is called {\it homogeneous} of degree $\deg f=h$ if
$f(V_g)\subseteq V_{hg}$ for all $g\in G$. Then the algebra ${\rm End}~ V$ 
of all linear transformations of $V$ becomes a $G$-graded algebra and
any matrix algebra $M_n(F)$ with an elementary $G$-grading is isomorphic to the endomorphism algebra ${\rm End}~ V$ of some $G$-graded vector space
(see, for example, \cite{BZ}, Proposition 2.2). In fact, any such 
isomorphism is given by fixing a homogeneous basis of $V$. Now it is
clear that the permutation of elements of a homogeneous basis of $V$ gives rise
to a graded isomorphism $R({\bf g})\rightarrow R({\bf g}_\sigma)$.

Throughout the paper, keeping in mind this isomorphism, we will identify 
the $G$-gradings on $M_n(F)$ defined by the $n$-tuples ${\bf g}=(g_1,\ldots,g_n)$
and ${\bf g}_\sigma=(g_{\sigma(1)},\ldots,g_{\sigma(n)})$.

\vskip .2in

\subsection{Tensor product of graded algebras}\label{ss2}
 Our second example allows one to construct new graded algebras from the given ones. First, let $G$ and $H$ be two 
groups while $A=\oplus_{g\in G}A_g$ and $B=\oplus_{h\in H}B_h$ the $G$-
and $H$-graded algebras, respectively. Then the tensor product 
$C=A\otimes B$ can be endowed with a $G\times H$-grading in the following
natural way. Given $t=gh\in G\times H$, we may set
$$
C_t=A_g\otimes B_h.
$$
Obviously, it is a well-defined $G\times H$-grading on $C$. A generalization
of this approach is as follows. Consider two $G$-graded algebras $A=\oplus_{g\in G}A_g$ and
$B=\oplus_{g\in G}B_g$ and set $S={\rm Supp}~A, T={\rm Supp}~B$. If the
elements of $S$ and $T$ pairwise commute then we can define a $G$-grading on
$C=A\otimes B$ by setting
$$
C_r=\sum_{gh=r}A_g\otimes B_h
$$
for all $r\in ST$. In particular, this definition is correct if $G$ is 
an abelian group. We cannot extend this definition to the non-abelian case in
general but there is a particular case when this can be done, and this definition is
very important for our classification.

\vskip .2in

\subsection{Induced grading on the tensor product}\label{ss3}
 Let $A=M_n(F)$ be a matrix algebra with an elementary $G$-grading defined by an $n$-tuple ${\bf g}=
(g_1,\ldots,g_n)$ and $B=\oplus_{g\in G}B_g$ any $G$-graded
algebra. We recall the notion of an {\it induced} grading on $C=A\oplus B$,
introduced in \cite{BZ}. Note that the same construction was studied in
\cite[Section 1.5]{NVO}, in the case where $B$ is a graded division algebra.

As in the commutative case, any homogeneous component of $C$ is defined as the span of the tensor
products of homogeneous elements but their degrees are defined differently, namely,
$$
C_g=\mathrm{Span}\,\{E_{ij}\otimes b\vert\deg b=h,~g_i^{-1}hg_j=g\}.
$$
Direct computations show that the decomposition $C=\oplus_{g\in G}C_g$ is
a $G$-grading and $B$ is a graded subalgebra of $C$. Also $A$ is a graded
subalgebra of $C$, provided that $B$ is a unitary algebra. Obviously,
$A\otimes B$ is isomorphic to the matrix algebra $M_n(B)$ if we identify
$E_{ij}\otimes b$ with the matrix $E_{ij}(b)$ with  only one non-zero
entry $b\in B$ on the intersection of the $i^{\mathrm{th}}$ row and the $j^{\mathrm{th}}$ column. Moreover, 
this is an isomorphism of graded algebras, if we define a $G$-grading on
$M_n(B)$ by setting
$$
\deg E_{ij}(b)= g_i^{-1}hg_j=g
$$
for any homogeneous $b\in B_h$.

Any grading on a matrix algebra produces a graded simple algebra. In the 
next subsection we present a series of graded algebras, which turn out to be graded
{\it division} algebras.

\vskip .2in

\subsection{Twisted group algebras}\label{ss4}
 Let $R=F[G]$ be a group algebra of
$G$ over a field $F$ that is $R$ is spanned by its basic elements 
$r_g, g\in G$, with the product $r_gr_h=r_{gh}$. Then $R$ is endowed
with the {\it canonical} $G$-grading $R=\oplus_{g\in G}R_g$ where 
$R_g=\mathrm{Span}\,\{r_g\}$ is a 1-dimensional vector space and all homogeneous
non-zero elements are invertible. Hence $F[G]$ is a graded division algebra.

Suppose now that the same vector space $\mathrm{Span}\,\{r_g\vert g\in G\}$ is endowed
with a different product:
\begin{equation}\label{r0}
r_gr_h=\sigma(g,h)r_{gh},
\end{equation}
where $\sigma(g,h)\in F^\ast$ is a non-zero scalar for any $g,h\in G$.
This product is completely determined by the mapping $\sigma: G\times G
\rightarrow F^\ast$. To be associative, this new product
should satisfy
\begin{equation}\label{r1}
\sigma(x,y)\sigma(xy,z)=\sigma(y,z)\sigma(x,yz)\mbox{ for all }x,y,z\in G
\end{equation}
(see, for example, \cite[Chapter 1]{P1}, for all details). The mapping
$\sigma: G\times G \rightarrow F^\ast$ satisfying (\ref{r1}) is called
a 2-{\it cocycle} on $G$ with values in $F^\ast$, while the associative algebra
$$
F^\sigma[G]=\mathrm{Span}\, \{r_g\vert g\in G\}
$$
with the product (\ref{r0}) is called a {\it twisted group algebra}
determined by $\sigma$. If, say, $\sigma\equiv 1$ then $F^\sigma[G]$ is
an ordinary group algebra.

Twisted group algebras inherit some properties of ordinary group algebras. 
For example, if we set $\deg r_g=g$ then, obviously, $F^\sigma[G]$ becomes
a $G$-graded algebra. We will call this grading {\it canonical}. Another 
important property is Maschke's theorem (see, for example, \cite[Theorem 4.4]{P2}).

\begin{Proposition}\label{prop1}Let $G$ be a finite group such that the order $|G|$ is
invertible in $F$. Then $F^\sigma[G]$ is semisimple for any 2-cocycle
$\sigma$. \epf
\end{Proposition}

\section{Unitarity of graded simple algebras}\label{s3}

Until the end of the paper $G$ will stand for the grading group and
$e$ for its identity element. The main goal of this section is to prove that any
graded simple algebra over an arbitrary field $F$ has the unit. It is known 
that any semisimple Artinian ring has the unit element (\cite{H}, 
Section 1.4). Following the similar way we will prove the unitarity of 
a graded simple finite-dimensional algebra. For shortness we will call
a homogeneous idempotents of a graded algebra as {\it graded idempotents}.

\begin{Lemma}\label{L2}
Let $R=\oplus_{g\in G}R_g$ be a finite-dimensional graded simple algebra and  $I\subset R$ a minimal graded right ideal of $R$. Then $I=aR$ for 
some graded idempotent $a$.
\end{Lemma}

{\em Proof.} First, suppose $I^2=0$. Since $(RI)^2=R(IR)I\subset RI^2=0$,  the two-sided ideal $RI$ 
of $R$ is nilpotent. Obviously, $RI$ 
is a  graded ideal and so by the hypothesis of our lemma we have $RI=0$. But this
means that $I$ is a non-trivial two-sided graded ideal of $R$, a
contradiction. Hence $I^2\ne 0$ and there exists a homogeneous $x\ne 0$
such that $xI\ne 0$. By the minimality of $I$ we have $xI=I$. In particular,
$xa=x$ for some non-zero homogeneous $a\in I$.

Now the right annihilator $X$ of $x$ in $R$ is a graded right ideal and 
either $I\cap X=0$ or $I\cap X=I$. Since $xa=x\ne 0$ it follows that $I\cap X=0$. Finally,
$xa^2=xa=x$ implies $x(a^2-a)=0$. Hence $a^2=a$ is a graded idempotent in 
$I$. Since $aI$ is a non-zero graded right ideal contained in $I$, we 
conclude that $aI=I$. \hfill $\Box$

\begin{Lemma}\label{L3}
Let $I\subset R$ be an arbitrary non-zero graded right ideal of a graded
simple finite-dimensional algebra $R$. Then $I=bR$ for some graded
idempotent~$b$.
\end{Lemma}

{\em Proof.} Since $\dim R<\infty$ a right ideal $I$ contains a minimal
ideal. Then by Lemma \ref{L2} it contains a graded idempotent $a$. Now, 
given a graded idempotent $t\in I$, we denote by $A(t)$ the right annihilator 
of $t$ in $I$:
$$
A(t)=\{x\in I\vert tx=0\}.
$$
We claim that there is an idempotent $b\in A$ such that $A(b)=0$. To prove this, it is sufficient to check that if $A(a)\ne 0$ then one can find another graded idempotent
$t\in I$ with $\dim A(t)<\dim A(a)$. 

Since $a$ is homogeneous, the right annihilator $T$ of $a$ in $R$ is a
graded right ideal. Therefore $A(a)=T\cap I$ is also a graded right ideal.
As before, applying Lemma \ref{L2} one can find a homogeneous idempotent
$f\in A(a)$. Then $f^2=f$ and $af=0$. Denote
$$
t=a+f-fa.
$$
Then
$$
t^2=(a+f-fa)(a+f-fa)=a^2+fa+f^2-f^2a-fa^2=a+f-fa=t,
$$
that is, $t$ is an idempotent. Note that any homogeneous idempotent of a
graded algebra $R$ lies in the identity component $R_e$. Hence, $a,f\in R_e$
and $t$ is also graded idempotent.

Let us check now that $A(t)$ is a proper subspace of $A(a)$. First pick
$x\in A(t)$. Then
$$
0=tx=atx=a(a+f-fa)x=a^2x=ax
$$
i.e. $x$ annihilates $a$. On the other hand, $af=0$ but $tf=(a+f-fa)f=f^2=
f\ne 0$. Hence $A(t)\ne A(a)$.

Continuing this process we obtain a graded idempotent $b\in I$ such that
$A(b)=0$, i.e. $bx\ne 0$ for any $0\ne x\in I$. Since $b(x-bx)=0$ for any 
$x\in I$ we get $x=bx$, in particular, $bI=I$. Recall that $b\in I$ and $I$ 
is a graded right ideal. Therefore, $bR\subset I$ and we got the required
equality $I=bR$, where $b=b^2$, is a graded idempotent. \hfill $\Box$

Now we can prove the existence of the unit element in any finite-dimensional
graded simple algebra.

\begin{Theorem}\label{T1} 
Let $R$ be a finite-dimensional graded simple algebra over an arbitrary
field. Then $R$ is a unitary algebra.
\end{Theorem}

{\em Proof.} By Lemma \ref{L3} there exist a homogeneous idempotent $a\in R$
such that $R=aR$. Obviously, $ax=x$ for all $x\in R$. Now let us set $I=\{x-xa\vert x\in R\}$. Since $a\in R_e$, $I$ is 
a graded subspace. Clearly, $I$ is a left ideal and $Ia=0$. Hence $IR=IaR=0$
and $I$ is a graded two-sided ideal. By the hypothesis, then either $I=0$ or
$I=R$. If $I=0$, then $R^2=IR=0$,  a contradiction. Hence $I=0$, 
that is, $xa=x$ for any $x\in R$. It follows that $xa=ax=x$ for all $x\in R$,
i.e. $a=1$ and the proof of the theorem is complete.
\hfill $\Box$

From now on we will denote the unit of $R$ by $1$.

\section{Graded division algebras}\label{s4}

In this Section we describe finite-dimensional graded division algebras over 
an arbitrary algebraically closed field. Obviously, some of these results are known, and we list them with the proofs here for easier reference.

\begin{Lemma}\label{D1}
Let $R=\oplus_{g\in G}R_g$ be a finite-dimensional graded division algebra
over an algebraically closed field $F$. Then $H=\su{R}$ is a subgroup
of $G$ and $\dim R_h=1$ for any $h\in H$.
\end{Lemma}

{\em Proof.} Our claim about the support is obvious since no invertible element
can be a zero divisor, therefore $R_{gh}\supset R_gR_h\ne 0$ as soon as
$R_g\ne 0,R_h\ne 0$, that is, $H$ is a finite multiplicatively closed subset
of the group $G$.

Now let $R$ be a graded division algebra over $F$, $\dim R<\infty$. Then 
$R_e$ is a finite-dimensional division algebra over an algebraically closed
field $F$. Hence $R_e=F$. In particular, $\dim R_e=1$.

Now let $g\ne e$ and $x\in R_g$ nonzero. Then $x$ has an inverse $x^{-1}$.
Clearly, $x^{-1}$ is homogeneous and $x^{-1}\in R_{g^{-1}}$. Pick an
arbitrary $y\in R_g$. Then $x^{-1}y\in R_e$ and $y=\lambda x$, for $\lambda\in F$, so that 
$\dim R_g=1$, and the proof of our lemma is 
complete.
\hfill $\Box$

Now we  are ready to characterize graded division algebras as twisted group
rings.

\begin{Theorem}\label{T2}
Let $R=\oplus_{g\in G}R_g$ be a finite-dimensional $G$-graded algebra over
an algebraically closed field $F$. Then $R$ is a graded division algebra
if and only if $R$ is isomorphic to the twisted group algebra $F^\sigma[H]$
with the canonical $H$-grading where $H$ is a finite subgroup of $G$ and
$\sigma: H\times H\rightarrow F^\ast$ is a 2-cocycle on $H$
\end{Theorem}

{\em Proof.} Obviously, any twisted group algebra is a graded division 
algebra. Conversely, let $R$ be a graded division algebra. Then 
${\rm Supp}~R$ is a finite subgroup of $G$ by Lemma \ref{D1} and also 
$\dim R_h=1$ for any $h\in H$. Fix an arbitrary non-zero $r_h\in R_h$.
Then, given $g,h\in H$, one has
$$
r_gr_h=\sigma(g,h)r_{gh}
$$
with some non-zero scalar $\sigma(g,h)\in F$ since $R_gR_h\in R_{gh}$ and
dim $R_{gh}=1$. Since $R$ is an associative algebra, the scalars $\sigma(g,h)$
should satisfy
$$
\sigma(x,y) \sigma(xy,z) = \sigma(y,z) \sigma(x,yz)
$$
(see Introduction), i.e. $\sigma: H\times H\rightarrow F^\ast$ is a 
2-cocycle. Then $R\cong F^\sigma[H]$ and the proof of the 
theorem is complete. \hfill $\Box$

At the end of section we characterize graded division algebras in the class
of all graded simple algebras.

\begin{Lemma}\label{D2}
Let $C=\oplus_{g\in G}C_g$ be a finite-dimensional graded simple algebra
over an arbitrary field $F$. If $\dim C_e=1$ then $C$ is a graded
division algebra.
\end{Lemma}

{\em Proof.} By Theorem \ref{T1}, $C$ is a unitary algebra. Clearly,
$1\in C_e$ and $C_e=F$. Let $x\in C_g,g\ne e$, be a non-zero homogeneous
element. Since $C$ is graded simple we have $CxC=C$. In particular, there
exist homogeneous $a,b\in C$ such that $0\ne axb\in C_e$. Without any loss
of generality, we may assume that $axb=1$. We recall that for arbitrary 
finite-dimensional unitary algebra over a field $F$ the equality
$uv=1$ implies $vu=1$. Hence $bax=xba=1$, i.e. $ba=x^{-1}$ and 
the proof is complete.
\hfill $\Box$

If $F$ is algebraically closed then by Theorem \ref{T2} a finite-dimensional
graded simple algebra $C$ is a graded division algebra if and only if
$\dim C_e=1$.

\section{Semiprimitivity of graded division algebras}\label{s5}

In this section we prove the graded 
simple algebras are semiprimitive, as non-graded algebras. First we study certain subalgebras of such algebras. We start
with an easy remark about the identity component of a graded simple
algebra.

\begin{Lemma}\label{L1}
Let $R=\oplus_{g\in G}R_g$ be a finite-dimensional $G$-graded simple algebra
over an arbitrary field. If $R$ is graded simple then $R_e$ is semiprimitive.
\end{Lemma}

{\em Proof.} In the previous paper \cite{ZS} authors proved that for any
$G$-graded ring $R=\oplus_{g\in G}R_g$ with $|{\rm Supp}~R|<\infty$ the
identity component $R_e$ does not contain nilpotent ideals as soon as $R$
itself has no non-zero nilpotent ideals (see Lemma 3). In fact, it was 
proved that if $I$ is a nilpotent ideal of $R_e$ then $RIR$ is a two-sided 
nilpotent ideal of $R$. Obviously, $RIR$ is a graded ideal, hence the graded 
simplicity of $R$ implies $I=0$ and $R_e$ is semiprimitive.
\hfill $\Box$

By the previous Lemma, a graded simple algebra $R$ contains at least one
graded idempotent. We use these idempotents for studying the structure of
$R$ and its subalgebras.

\begin{Lemma}\label{S1}
Let $R=\oplus_{g\in G}R_g$ be a finite-dimensional graded simple algebra
and $t=t^2\in R_e$ a graded idempotent. Then $tRt$ is a graded simple algebra.
\end{Lemma}

{\em Proof.} Denote $A=tRt$. By Theorem \ref{T1}, $R$ is a unitary algebra.
The statement of Lemma is obvious if $t=1$. Suppose $t\ne 1$. Since $t\in R_e$, it is clear that $A$ is a graded subalgebra of $R$. Since $R$ is unitary, $t\in A$, so that $A$ is nonzero. We only need to prove that $A$ is graded simple.



Let $I\subset A$ be a graded ideal of $A$. We generate by $I$ a 
two-sided ideal $T=RIR$ of $R$. First note that $T\cap A=I$.

Indeed, if $y\in T\cap A$ then
$$
y=\sum_i r_ix_is_i\in A
$$
where $x_i\in I$, $r_i, s_i\in R$ for all possible values of $i$.
Since $a=tat$ for any $a\in I$, $tr_it,ts_it$ are in $A$ and $I$ is an ideal of $A$ we have
$$
y=tyt=\sum (tr_it)x_i(ts_it)\in I.
$$
 That is 
$T\cap A\subset I$. The inverse containment is obvious.

It follows that any proper graded ideal of $A$ generates a proper graded
ideal of $R$. We refer to the simplicity of $R$ and thus  the proof is complete.\epf

Recall that by Lemma \ref{L1} the identity component $R_e$ of a graded simple
algebra $R$ is semiprimitive, hence the sum of simple ideals. To prove that $R$ is itself semiprimitive, we first
consider the case when $R_e$ is simple.

\begin{Lemma}\label{S2}
Let $R=\oplus_{g\in G}R_g$ be a finite-dimensional graded simple algebra
over an algebraically closed field $F$ such that $A=R_e$ is simple. Then 
$R=AC\cong A\otimes C$ where $A\cong M_k(F)$ with the trivial $G$-grading, $C=\oplus_{g\in G}C_g$ is the centralizer of $A$ in $R$ and $C$ is a graded division algebra.
\end{Lemma}

{\em Proof.} By Theorem \ref{T1}, $R$ is unitary. Since $F$ is
algebraically closed, $A$ is isomorphic to $M_k(F)$ and $1\in A$. Denote by $C$ the centralizer of $A$ in $R$. Since 
$A=R_e$, $C$ is a graded subalgebra. Moreover, $A\cap C=F$. By \cite[Lemma 3.11]{S} or \cite[Chapter 4, Section 3]{J}, one has $R=AC\cong A\otimes C$. Since $A\cap C=F$ we also 
have $\dim C_e=1$. Then $C$ is a graded division algebra by Lemma \ref{D2}.\epf

Now we have all ingredients for proving the semiprimitivity of a graded simple 
algebra. Note that if $G$ is finite and $|G|$ is invertible in $F$ then
the Jacobson radical of any $G$-graded $F$-algebra with $1$ is a graded ideal 
(see \cite[Theorem 4.4]{CM}).In particular, such graded simple algebra are semiprimitive. We need to extend this latter statement to 
an arbitrary group $G$.

\begin{Lemma}\label{S3}
Let $R=\oplus_{g\in G}R_g$ be a finite-dimensional graded simple algebra
over an algebraically closed field $F$ such that ${\rm char}~F=0$ or
${\rm char}~F$ is coprime with the order of any finite subgroup of $G$. Then 
$R$ is semiprimitive.
\end{Lemma}

{\em Proof.} By Lemma \ref{L1}, $ R_e$ is semiprimitive. Then $R_e=A^{(1)} 
\oplus\ldots\oplus A^{(m)}$ where all summands $A^{(i)}$ are matrix
algebras. Let $e_1,\ldots, e_m$ be the units of $A^{(1)}$,\ldots, $A^{(m)}$,
respectively. By Theorem \ref{T1} $R$ is a unitary algebra, hence
$e_1+\ldots+e_m=1$. Clearly $e_1,\ldots, e_m$ is the complete system of orthogonal idempotents,
hence
$$
R=\oplus_{1\le i,j\le m} e_iRe_j.
$$
First we fix one component $B=e_iRe_i$. Then, by Lemma \ref{S1}, $B$  is graded
simple and then, by Lemma \ref{S2}, $B\cong M_k(F)\otimes C$  with a graded
division algebra $C$. By Theorem \ref{T2} $C$ is isomorphic to a 
twisted group algebra $F^\sigma[H]$ for a finite subgroup $H\subset G$
and a 2-cocycle $\sigma: H\times H\rightarrow F^\ast$. It is known that
$F^\sigma[H]$ is semiprimitive if $|H|^{-1}\in F$ (see, for example, \cite[Theorem 4.2]{P2}
). Since $B\cong M_n(F^\sigma[H])$, the Jacobson radical of $B$ is known \cite{H} to be $M_n(J)$ where $J$ is the radical of $F^\sigma[H]$. Therefore $B$ is also semiprimitive.

We have proved that all components $e_iRe_i$ are semiprimitive subalgebras.
Denote now by $J=J(R)$, the Jacobson radical of $R$. Then $J\cap e_iRe_i=e_iJe_i=0$ for
all $i=1,\ldots,m$. Suppose now that $J$ is non-zero. Since $J=\oplus_{i,j}e_iJe_j$, it follows that $e_kJe_l\ne 0$ for
some $k\ne l$. Fix some 
$0\ne x\in e_kJe_l$. Since
$$
x=e_kxe_l=\sum_{g\in G} e_kx_ge_l
$$
and since $e_kx_ge_l$ is homogeneous in the $G$-grading, with $\deg(e_kx_ge_l)=g$,
it follows that $x_g=e_kx_ge_l$, that is, all homogeneous components of $x$ are in
$e_kJe_l$. In particular, $e_ix_ge_j=0$ as soon as either $i\ne k$ or $j\ne l$.

Now let us consider the set $T=e_iRe_k x e_lRe_i$. Clearly, $T\subset e_iRe_i 
\cap J$, hence $T=0$. But then $e_iRe_k x_g e_lRe_i=0$ for any component 
$x_g$. Finally we obtain
$$
Rx_gR= \bigoplus_{i,j,s,t} e_iRe_j x_g e_sRe_t =
\bigoplus_{i\ne j} e_iRe_k x_g e_lRe_j,
$$
that is, $x_g$ generates in $R$ a proper graded ideal, a contradiction.
Hence $J=J(R)=0$ and the proof of the lemma  is now complete.\epf

\section{Description of graded simple algebras}\label{s6}

First we need one technical remark from Linear Algebra.

\begin{Lemma}\label{G1}
Let $B=B_1\oplus\ldots\oplus B_n$ be the direct sum of matrix algebras and $z_1,z_2\in B$ two orthogonal idempotents. Consider the decomposition
$z_i=z_i^1+\ldots+z_i^n$, $i=1,2$, with $z_i^j\in B_j$. Then
$$
\dim z_1Bz_2= \sum_{i=1}^n ({\rm rank}~ z_1^i)({\rm rank}~ z_2^i)
$$
where ${\rm rank}~ z_i^j$ is an ordinary matrix rank in $B_j$.
\end{Lemma}

{\em Proof.} Clearly, any $z_i^j$ is an idempotent of $B_j$ and
$z_1^j z_2^j = z_2^j z_1^j=0$ since $z_1$ and $z_2$ are orthogonal. Then,
up to a conjugation by an inner automorphism of $B_j$, these elements are
$z_1^j=E_{11}+\ldots+E_{mm}, z_2^j=E_{m+1,m+1}+\ldots+E_{m+k,m+k}$ where
$E_{ii}$'s are diagonal matrix units of $B_j$ and ${\rm rank}~ z_1^j=m,
{\rm rank}~z_2^j=k$. Then, obviously,
$$
\dim z_1^jB_jz_2^j= mk=({\rm rank}~ z_1^j)({\rm rank}~ z_2^j).
$$
Now our lemma follows from the equality $z_1Bz_2=z_1^1B_1z_2^1\oplus \ldots
\oplus z_1^nB_nz_2^n$.

\hfill $\Box$

Now we consider a graded simple algebra with the identity component consisting
of two simple summands.

\begin{Lemma}\label{G2}
Let $R=\oplus_{g\in G}R_g$ be a finite-dimensional graded simple algebra
over an algebraically closed field $F$ such that either ${\rm char}~F=0$ or
${\rm char}~F$ is coprime to the order of any finite subgroup of $G$.
Let the identity component of $R$ be the sum of two simple summands, 
$R_e=A_1\oplus A_2$. Let $e_1\in A_1$ and $e_2\in A_2$
be the units of $A_1, A_2$ respectively. Set $R_1=e_1Re_1, 
R_2=e_2Re_2$ and consider the decompositions $R_1=A_1C_1, 
R_2=A_2C_2$ given in Lemma \ref{S1}, with graded division algebras
$C_1,C_2$. Let $d_1\in A_1$ and $d_2\in A_2$ be a pair of minimal
idempotents such that $M=d_1Rd_2\ne 0$. Then
\begin{itemize}
\item $\dim C_1=\dim C_2=\dim M$;

\item $M=C_1x=xC_2$ for any homogeneous $0\ne x\in M$.
\end{itemize}
\end{Lemma}

{\em Proof.} By Lemma \ref{S3}, $R$ is semiprimitive. Let 
$R=B_1\oplus\cdots\oplus B_n$ be the decomposition of $R$ as the sum of 
simple ideals, as a non-graded algebra. According to this decomposition, we
write the unit of $A_1$ as $e_1=e_1^1+\ldots+e_1^n, e_1^i\in B_i$,
$i=1,\ldots, n$. Then
\begin{equation}\label{onestar}
e_1Re_1=e_1^1 B_1 e_1^1\oplus\cdots\oplus e_1^n B_n e_1^n=A_1C_1.
\end{equation}
Since $C_1$ commutes with $A_1$  we have 
$C_1=C_1^1\oplus\cdots\oplus C_1^n$ where $C_1^i= C_1\cap B_i$. Then we fix
$1\le i\le n$ and denote by $\varphi$ the canonical projection of $B$ onto
$B_i$. Let $A_1^i=\varphi(A_1)$. Since $A_1$ is simple, then either $A_1^i=0$ or $A_1^i\cong A_1$. If $A_1^i=0$ then 
$$
A_1\subset B_1\oplus\cdots\oplus B_{i-1}\oplus B_{i+1}
\oplus\cdots\oplus B_n=R'
$$
and $e_1Re_1$ is a graded subspace of $R$ contained in a proper 
non-graded ideal of $R$. Hence $e_1Re_1$ generates a non-trivial graded
ideal of $R$, a contradiction. It follows that $A_1^i\cong A_1$ is simple 
and $C_1^i$ commutes with $A_1^i$.

Denote $B_i'=e_1^i B_i e_1^i$. Since $e_1^i=
\varphi(e_1)$ is an idempotent, $B_i'$ is simple . Also $C_1^i\subset B_i'$. Then it follows from (\ref{onestar})  that
$B_i'=A_1^iC_1^i$ and $C_1^i$ is the centralizer of $A_1^i$ in $B_i'$. Hence
$C_1^i$ is simple, $C_1^i\cong M_{p_i}(F)$ for some $p_i\ge 1$ and
\begin{equation}\label{plus}
B_i'\cong A_1^i\otimes C_1^i.
\end{equation}

If we decompose $d_1=d_1^1+\ldots+d_1^n$ where $d_1^1\in B_1,\ldots,
d_1^n\in B_n$ then $d_1^i=\varphi(d_1)$ is a minimal idempotent of $A_1^i$
and from (\ref{plus}) it follows that the rank of $d_1^i$ in $B_i'$ is
equal to $p_i$. Note that $B_i'\subset B_i$ are two matrix algebras and
$B_i'=e_1^iB_ie_1^i$ where $e_1^i$ is an idempotent. Hence for any 
$x\in B_i'$ its rank as a matrix in $B_i'$ coincides with its rank
in $B_i$. It follows that ${\rm rank}~d_1^i=p_i$.

Similarly, $C_2=C_2^1\oplus\ldots\oplus C_2^n$, $C_2^j\cong M_{q_j}(F)$,
$1\le j \le n$, $d_2=d_2^1+\ldots+d_2^n$, $d_2^j\in B_j$, $1\le j\le n$,
and ${\rm rank}~d_2^i=q_i$. Then by Lemma \ref{G1}
\begin{equation}\label{twostars}
\dim M= \sum_{i=1}^n ({\rm rank}~ d_1^i)({\rm rank}~ d_2^i)
\end{equation}
and
$$
\dim C_1=p_1^2+\ldots+p_n^2,\quad \dim C_2=q_1^2+\ldots+q_n^2.
$$

Since $C_1$ is a graded simple division algebra, the dimension of any non-trivial graded left $C_1$-module
is at least $\dim C_1$, and since $d_1$ centralizes $C_1$ we know that
$M=d_1Rd_2$ is a left $C_1$-module.
Hence $\dim M\ge \dim C_1$. Similarly, $\dim M\ge \dim C_2$. Comparing with
(\ref{twostars}) we obtain
$$
p_1q_1+\ldots+p_nq_n \ge p_1^2+\ldots +p_n^2,
$$
$$
p_1q_1+\ldots+p_nq_n \ge q_1^2+\ldots +q_n^2.
$$
Hence
$$
2\sum p_iq_i \ge \sum(p_i^2+q_i^2) \quad {\rm and}\quad 
\sum(p_i-q_i)^2\le 0.
$$

It follows that $p_1=q_1,\ldots,p_n=q_n$, i.e. $\dim C_1=\dim C_2$ and
then  by (\ref{twostars}) $\dim M=\dim C_1=\dim C_2$.

Finally, if we take any homogeneous $x\ne 0$ in $M$ then $C_1x$ is a nonzero graded
subspace of $M$. Since $C_1$ is a graded division
algebra, we must have $\dim M=\dim C_1$ . Hence $C_1x=M$. Similarly, $xC_2=M$ and the proof is complete.\epf

Now we are in a position to prove the main result of our paper.

\begin{Theorem}\label{T3}
Let $R=\bigoplus_{g\in G}R_g$ be a finite-dimensional algebra over an 
algebraically closed field $F$ graded by a group $G$. Assume that either 
${\rm char}~F=0$ or ${\rm char}~F$ is coprime to the order of any finite 
subgroup of $G$. Then $R$ is graded simple if and only if $R$ is isomorphic 
to $M_k(F)\otimes F^\sigma[H]\cong M_k(F^\sigma[H])$, a matrix algebra
over the graded division algebra $F^\sigma[H]$ where $H$ is a finite subgroup of $G$ and $\sigma: H\times H\rightarrow F^\ast$ is a 2-cocycle on 
$H$. The $H$-grading on $F^\sigma[H]$ is canonical and the $G$-grading on
$M_k(F^\sigma[H])$ is defined by a $k$-tuple $(a_1,\ldots,a_k)\in G^k$, so that $\deg(E_{ij}\otimes x_h)= a_i^{-1}ha_j$ for any matrix unit 
$E_{ij}$ and any homogeneous element $x_h\in F^\sigma[H]_h$.
\end{Theorem}

{\em Proof.} The proof that $M_k(F^\sigma[H])$ is graded simple algebra is
easy and left to the reader. The main part is the proof that any graded simple
algebras is a full matrix algebra over a graded division algebra.

So, we assume that $R$ is graded simple. Then $R$ is semiprimitive by Lemma \ref{S3},
$A=R_e$ is semiprimitive by Lemma \ref{L1} and we can write
$A=A^{(1)}\oplus\ldots\oplus A^{(m)}$ with $A^{(i)}\cong M_{n_i}(F)$, 
$i=1,\ldots,m$. Denote by $e_1,\ldots,e_m$ the units of $A^{(1)},\ldots,
A^{(m)}$, respectively. Then $e_1,\ldots,e_m$ are orthogonal idempotents and 
$1=e_1+\ldots+e_m$. First we will find a graded subalgebra of $R$ isomorphic 
to $M_m(F)$ with an elementary $G$-grading.

Denote $R^{(i)}=e_iRe_i$. Then $R^{(i)}$ is a graded subalgebra of $R$.
Moreover by Lemma \ref{S1} $R^{(i)}$ is graded simple.  By Lemma \ref{S2} $R^{(i)}=
A^{(i)}C^{(i)}\cong A^{(i)}\otimes C^{(i)}$  where 
$C^{(i)}$ is the centralizer of $A^{(i)}$ in $R^{(i)}$ and $C^{(i)}$ is a
graded division algebra. We first show that there exist homogeneous
$x_{12}\in e_1Re_2,\ldots,x_{m-1,m}\in e_{m-1}Re_m$ and
$x_{21}\in e_2Re_1,\ldots,x_{m,m-1}\in e_{m}Re_{m-1}$ such that
\begin{equation}\label{*1}
w_k=e_kx_{k,k-1}e_{k-1}\cdots e_2 x_{21}e_1x_{12}e_2\cdots
x_{k-1,k}e_k\ne 0
\end{equation}
for all $k=2,\ldots,m$ and $\deg w_k=e$ in $G$-grading.

Note that
\begin{equation}\label{*2}
W_k=e_kR\cdots Re_2Re_1Re_2R\cdots Re_k = R^{(k)}\ne 0
\end{equation}
for all $k=2,\ldots,m$. This follows by induction on $k$ with obvious basis for $k=1$, where $W_1=\{ e_1\}$. Now if $W_{k-1}=R^{(k-1)}\ne 0$ then the ideal $RW_{k-1}R$ is graded in $R$, so that $RW_{k-1}R=R$. Therefore, $e_kRW_{k-1}Re_k = e_kRe_k=R^{(k)}\ne 0$, as expected.

We will prove (\ref{*1}) by induction on $k$. If $k=2$ then by 
(\ref{*2}) there exist homogeneous $x_{21}\in e_2Re_1, x_{12}\in e_1Re_2$
such that $w_2=e_2x_{21}e_1x_{12}e_2$ is non-zero. Denote $H^{(2)} = 
{\rm Supp}~C^{(2)}$. Then by Lemma \ref{D1} $H^{(2)}$ is a subgroup of $G$.
Moreover, ${\rm Supp}~R^{(2)}=H^{(2)}$ following since $R^{(2)}=A^{(2)}C^{(2)}$ and
$A^{(2)}\subset R_e$. Now, if $\deg w_2=h\ne e$ then there exist
$0\ne b\in C^{(2)}$ with $\deg b=h^{-1}$ and we can replace $x_{21}$ by 
$x_{21}'=bx_{21}$. Then $w_2'=bw_2=e_2bx_{21}e_1x_{12}e_2=
e_2x_{21}'e_1x_{12}e_2\in R_e$ is non-zero since $b$ is invertible.
Similarly, if $w_{k-1}$ has been found then there exist homogeneous
$x_{k,k-1},x_{k-1,k}$ such that $w_k=e_kx_{k,k-1}w_{k-1}x_{k-1,k}e_k\ne 0$
and by the same argument as before one can take $x_{k,k-1}$ such that 
$\deg w_k=e$.

In particular we get
\begin{equation}\label{*3}
e_mx_{m,m-1}e_{m-1}\cdots e_2 x_{21}e_1x_{12}e_2\cdots
x_{m-1,m}e_m\ne 0
\end{equation}
and $\deg e_kx_{k-1,k}\cdots x_{k-1,k}e_k=e$ in $G$ for all $k=2,\ldots,m$.

We now denote by $E_{\alpha,\beta}^i$ the matrix units of $A^{(i)}$. Then by (\ref{*3})
there exist some integers $\alpha_1,\ldots,\alpha_m, \beta_1,\ldots,\beta_m$
such that
$$
E^m_{\alpha_m,\alpha_m} x_{m,m-1} E^{m-1}_{\alpha_{m-1},\alpha_{m-1}}
\cdots E^2_{\alpha_2,\alpha_2} x_{21} E^1_{\alpha_1,\beta_1} x_{12}
E^2_{\beta_2,\beta_2}\cdots 
$$
\begin{equation}\label{*4}
\cdots E^{m-1}_{\beta_{m-1},\beta_{m-1}} x_{m-1,m}
E^m_{\beta_m,\beta_m} \ne 0.
\end{equation}
Actually, in the above equation $\alpha_1=\beta_1$.

Now we define $y_{i,i+1}$ and $y_{i+1,i}$ by setting
$$
y_{i,i+1}=E^i_{1,\beta_i}E^i_{\beta_i,\beta_i} x_{i,i+1}
E^{i+1}_{\beta_{i+1},\beta_{i+1}}E^{i+1}_{\beta_{i+1},1},
\;y_{i+1,i}=E^{i+1}_{1,\alpha_{i+1}}E^{i+1}_{\alpha_{i+1},\alpha_{i+1}} 
x_{i+1,i}E^{i}_{\alpha_{i},\alpha_{i}}E^{i}_{\alpha_{i},1}
$$
for all $i=1,\ldots,m-1$. Then all $y_{i,i+1},y_{i+1,i}$ are homogeneous
and
\begin{equation}\label{*5}
E_{11}^iy_{i,i+1}E_{11}^{i+1}=y_{i,i+1},\quad
E_{11}^{i+1}y_{i+1,i}E_{11}^{i}=y_{i+1,i}.
\end{equation}
Moreover, by (\ref{*4}) we have
\begin{equation}\label{*6}
y_{m,m-1}\cdots y_{21}y_{12}\cdots y_{m-1,m}\ne 0.
\end{equation}

Consider an element $y_{22}=y_{21}y_{12}\in A^{(2)}$. It is non-zero by (\ref{*6}) and $E_{11}^2y_{22}E_{11}^2=y_{22}$ in $A^{(2)}$. Hence
$y_{22}=\lambda E_{11}^2$. If we replace $y_{21}$ with $\lambda^{-1}y_{21}$
then (\ref{*5}) and (\ref{*6}) still hold and $y_{22}=E_{11}^2$. Now let
$y_{11}=y_{12}y_{21}$. Then $y_{11}^3=y_{11}^2$ since $y_{22}$ is an 
idempotent. As before, $y_{11}=\lambda E_{11}^1$ and hence $\lambda=1$ or 
$0$. Let us check that $y_{11}\ne 0$.

Suppose $y_{11}=y_{12}y_{21}=0$. By Lemma 6 the
       subalgebra $(e_1+e_2)R(e_1+e_2)$ of $R$ is graded simple. If we denote $d_2=E_{11}^2, d_1=E_{11}^1$
then $d_2Rd_1=d_2(e_1+e_2)R(e_1+e_2)d_1=y_{21}C^{(1)}$ by Lemma \ref{G2} hence $y_{12}Rd_1=
y_{12}d_2Rd_1=0$. If we decompose $e_1$ to the sum of minimal idempotents
$E_{11}^1,\ldots,E_{n_1,n_1}^1$ with $d_1=E_{11}^1$ then 
$E_{11}^1,\ldots,E_{n_1,n_1}^1$, $e_2,\ldots,e_m$ is a system of orthogonal idempotents whose sum is $1$. Therefore
\begin{equation}\label{*7}
Ry_{12}R\cap d_1Rd_1=d_1Ry_{12}Rd_1.
\end{equation}
On the other hand, the right hand side of (\ref{*7}) is zero since 
$y_{12}Rd_1=0$. Hence $Ry_{12}R$ is a proper graded ideal of $R$, a 
contradiction.

We have proved that $y_{11}=E_{11}^1$. Similarly, replacing $y_{32},\ldots,
y_{m,m-1}$ by corresponding scalar multiples we may assume that
\begin{equation}\label{*8}
y_{i,i+1}y_{i+1,i}= E_{11}^i,\quad y_{i+1,i}y_{i,i+1}= E_{11}^{i+1}
\end{equation}
for all $i=1,\ldots,m-1$. Now for any $1\le i<j\le m$ we set
$$
y_{ij}=y_{i,i+1}\cdots y_{j-1,j}, ~~ 
y_{ji}=y_{j,j-1}\cdots y_{i+1,i}.
$$
Then all $y_{ij}, 1\le i,j\le m$ are non-zero by (\ref{*6}) and from
(\ref{*5}) it easily follows that they are linearly independent. Using 
(\ref{*8}) and the definition of $y_{ij}$ it is easy to show that
$$
y_{ij}y_{rt}=\delta_{jr}y_{it}
$$
where $\delta_{jr}$ is Kronecker delta. Therefore the linear span of
$\{y_{ij}|1\le i,j\le m\}$ is a graded subalgebra of $R$ isomorphic to 
$M_m(F)$.

Given $i\ne j$, consider an idempotent $t=e_i+e_j$. By Lemma \ref{S1}
 $R'=tRt$ of $R$ is a graded simple subalgebra and $R_e'=A^{(i)}\oplus
A^{(j)}$. Applying Lemma \ref{G2} to $R'$ we see that
\begin{equation}\label{*9}
C^{(i)}y_{ij}=y_{ij}C^{(j)}.
\end{equation}
Since all homogeneous components of $C^{(i)}$ and $C^{(j)}$ are one-dimensional, (\ref{*9}) gives a well-defined mapping $\nu:
C^{(i)}\rightarrow C^{(j)}$ which is an isomorphism of non-graded algebras. 
Moreover, for $H^{(i)}={\rm Supp}~C^{(i)}, H^{(j)}={\rm Supp}~C^{(j)}$ we 
obtain that the subgroups $H^{(i)}$ and $H^{(j)}$ are conjugate in $G$ and 
$\nu$ is in fact an isomorphism of graded algebras, if we identify $H^{(i)}$
and $H^{(j)}$ up to this conjugation. Denote by $\varphi_i:
C^{(1)}\rightarrow C^{(i)}$ the corresponding isomorphism for $i=2,\ldots,m$
and set $\varphi_1=\mathrm{id}_{C^{(1)}}$ on $C^{(1)}$. Then
$$
\varphi_1(x)y_{1i}=y_{1i}\varphi_i(x)
$$
for all $i=1,\ldots,m$. Using this relation we will prove that
\begin{equation}\label{*10}
\varphi_i(x)y_{ij}=y_{ij}\varphi_j(x)
\end{equation}
for all possible $i,j$.

First, we assume $i<j$. Then $xy_{1i}=y_{1i}\varphi_i(x)$ hence 
$xy_{1j}$=$xy_{1i}y_{ij}$=$y_{1i}\varphi_i(x)y_{ij}$. On the other hand,
$xy_{1j}$=$y_{1j}\varphi_j(x)$=$y_{1i}y_{ij}\varphi_j(x)$. Multiplying on the
left by $y_{i1}$ and recalling that $y_{ii}\in A^{(i)}$ commutes with 
$\varphi_i(x)\in C^{(i)}$ we get
\begin{eqnarray*}
y_{i1}y_{1i}\varphi_i(x)y_{ij}=y_{ii}\varphi_i(x)y_{ij}=
\varphi_i(x)y_{ii}y_{ij}&=&\varphi_i(x)y_{ij},\\
y_{i1}y_{1i}y_{ij}\varphi_j(x)=y_{ii}y_{ij}\varphi_j(x)&=&
y_{ij}\varphi_j(x),
\end{eqnarray*}
that is,  (\ref{*10}) holds for any $i<j$. If $i>j$ then for any $k\ge 2$
$$
y_{11}x=xy_{11}=xy_{1k}y_{k1}=y_{1k}\varphi_k(x)y_{k1}
$$
and multiplying by $y_{k1}$ on the left we get
$$
y_{k1}x=y_{kk}\varphi_k(x)=\varphi_k(x)y_{k1}.
$$
Now the same argument as before give us (\ref{*10}) also for $i>j$.

Now we extend $\mathrm{Span}\,\{y_{ij}|1\le i,j\le m\}$ to a larger matrix 
subalgebra of $R$ with an elementary $G$-grading and also construct a non-graded subalgebra $C$ isomorphic to $C^{(1)}$ in the sense of non-graded algebras.

Given $x\in C^{(1)}$, we denote
$$
\bar x=x+\varphi_2(x)+\ldots+\varphi_m(x).
$$
Since $\varphi_i(x)\in e_iRe_i$ and $e_1,\ldots,e_m$ are orthogonal,
$C=\mathrm{Span}\,\{\bar x|x\in C^{(1)}\}$ is isomorphic to $C^{(1)}$. Moreover,
\begin{equation}\label{*11}
\bar x z_i=\varphi_i(x)z_i, z_i\bar x=z_i\varphi_i(x)
\end{equation}
for any $z_i\in e_iRe_i$. Recall that $A^{(i)}\cong M_{n_i}$ and 
$E_{\alpha,\beta}^i$'s are the matrix units of $A^{(i)}$. It is easy to see that
the elements
\begin{equation}\label{*11+}
E_{\alpha,1}^iy_{ij}E_{1\beta}^j, 1\le\alpha\le n_i, 1\le\beta\le n_j,
i,j=1,\ldots, m
\end{equation}
are linearly independent and homogeneous. Set $k=n_1+\ldots+n_m$ and consider a linear map
$\varphi: B\rightarrow M_k(F)$ be  from 
$B=\mathrm{Span}\,\{E_{\alpha,1}^iy_{ij}E_{1\beta}^j\}$ onto $M_k(F)$ defined on the
basis elements by
\begin{equation}\label{*12}
\varphi(E_{\alpha,1}^iy_{ij}E_{1\beta}^j)= E_{\mu\nu}
\end{equation}
with $\mu=n_1+\ldots+n_{i-1}+\alpha, \nu=n_1+\ldots+n_{j-1}+\beta$. Then
$\varphi$ is an algebra isomorphism, i.e. $B$ is a graded subalgebra of $R$
isomorphic to $M_k(F)$.

On the other hand, it follows from (\ref{*10}) and (\ref{*11})  that
$$
\bar x E_{\alpha,1}^iy_{ij}E_{1\beta}^j=
\varphi_i(x)E_{\alpha,1}^iy_{ij}E_{1\beta}^j=
E_{\alpha,1}^i\varphi_i(x)y_{ij}E_{1\beta}^j=
$$
$$
E_{\alpha,1}^iy_{ij}\varphi_j(x)E_{1\beta}^j=
E_{\alpha,1}^iy_{ij}E_{1\beta}^j\varphi_j(x)=
E_{\alpha,1}^iy_{ij}E_{1\beta}^j\bar x
$$
for any $\bar x\in C$ meaning that  $C$ centralizes $B$ in $R$. In particular, 
$C\cap B=F$ and hence $BC\cong B\otimes C$ as a non-graded algebra.

We will prove that $BC=R$. Fix some $1\le i,j\le m$ and set
$$
P_{ij\alpha\beta}=\mathrm{Span}\,\{ E_{\alpha,1}^iC^{(i)}y_{ij}E_{1\beta}^j\}.
$$
Note that $P_{ij\alpha\beta}=C^{(i)}P_{ij\alpha\beta}= 
P_{ij\alpha\beta}C^{(j)}$. Therefore
$$
Q_{ij}=\sum_{1\le\alpha\le n_i\atop{1\le\beta\le n_j}} P_{ij\alpha\beta}
\subset e_iRe_j
$$
is a left $A^{(i)}C^{(i)}$-module and right $A^{(j)}C^{(j)}$-module. Let us prove that $Q_{ij}=e_iRe_j$. If not, $Q_{ij}$ generates a graded ideal of 
$R$ such that
\begin{eqnarray*}
RQ_{ij}R\cap e_iRe_j&=&(e_iRe_i)Q_{ij}(e_jRe_j)=R^{(i)}Q_{ij}R^{(j)}=
A^{(i)}C^{(i)}Q_{ij}A^{(j)}C^{(j)}\\&=&Q_{ij}\ne e_iRe_j,
\end{eqnarray*}
 a contradiction. It
follows that $Q_{ij}=e_iRe_j$ and $BC=R$.

Using this equality, we extend the mapping (\ref{*12}) to an isomorphism $\varphi:R\rightarrow M_k(C^{(1)})$ of graded algebras where $C^{(1)}$ is
a graded division algebra. We first describe the grading on $B$. A matrix 
subalgebra of $B$ spanned by all $y_{ij},1\le i,j\le m,$ is graded and all
matrix units $y_{ij}$ are homogeneous. Hence by \cite[Lemma 1]{ZS} its 
grading is elementary and there exist $g_1,\ldots, g_m$ such that 
$\deg y_{ij}=g_i^{-1}g_j$. Moreover,  since 
$(ag_1,\ldots, ag_m)$ defines the same grading for any $a\in G$, we can assume that $g_1=e$. 
Then the $G$-grading on $B$ is also elementary and 
$\deg(E_{\alpha,1}^iy_{ij}E_{1\beta}^j)=g_i^{-1}g_j$ for any suitable $\alpha,\beta$. In fact this $G$-grading on $B$ is defined by a $k$-tuple
$$
(a_1,\ldots,a_k)=(\underbrace{g_1,\ldots,g_1}_{n_1},\ldots,
\underbrace{g_m,\ldots,g_m}_{n_m}).
$$
Now we can define an elementary grading on $M_k(F)$ by setting $\deg E_{\mu\nu}=
a_\mu^{-1}a_\nu$. Then $\varphi: B\rightarrow M_k(F)$ given by (\ref{*12})
is an isomorphism of graded algebras. Finally we define a $G$-grading on 
$M_k(C^{(1)})$. We set $H=H^{(1)}={\rm Supp}~C^{(1)}$. We know that $H$ is a subgroup of $G$. As 
a non-graded algebra $M_k(C^{(1)})$ is isomorphic to $M_k(F)\otimes C^{(1)}$.
If we set $\deg(E_{\mu\nu}\otimes y_h)=a_\mu^{-1}ha_\nu$ for any homogeneous
$y_h\in C^{(1)}_h$ then direct computations show that this is a well-defined
grading. As it was remarked before this theorem, this is exactly the induced grading
on the tensor product on $M_k(F)\otimes C^{(1)}$.

In order to extend the isomorphism $\varphi$ from (\ref{*12}) to the map
$R=BC\rightarrow M_k(C^{(1)})$ we set
\begin{equation}\label{*13}
\varphi(\bar x_hE_{\alpha,1}^iy_{ij}E_{1\beta}^j)= E_{\mu\nu}\otimes x_h
\end{equation}
for any $\bar x_h\in C$ with $x_h\in C^{(1)}_h$ where 
$\mu=n_1+\ldots+n_{i-1}+\alpha, \nu=n_1+\ldots+n_{j-1}+\beta$. Clearly,
$\varphi$ is an isomorphism in a non-graded sense. We need to check that
$\varphi$ preserves grading.

The element on the right hand side of (\ref{*13}) is homogeneous and its
degree is $a_\mu^{-1}ha_\nu=g_i^{-1}hg_j$. Although $\bar x_h\in C$ is
a non-graded element of $R$, the product 
$$
w=\bar x_hE_{\alpha,1}^iy_{ij}E_{1\beta}^j=
E_{\alpha,1}^i\varphi_i(x_h)y_{ij}E_{1\beta}^j
$$
is homogeneous and  $\deg w=\deg(\varphi_i(x_h))\deg y_{ij}=
\deg(\varphi_i(x_h))g_i^{-1}g_j$. To compute this degree, we recall that
$x_hy_{1i}=y_{1i}\varphi_i(x_h)$ and $\deg y_{1i}=g_1^{-1}g_i=g_i$ since
$g_1=e$. Hence 
$$
\deg(x_hy_{1i})=hg_i=\deg(y_{1i}\varphi_i(x_h))=g_i\deg\varphi_i(x_h)
$$
and $\deg\varphi_i(x_h)=g_i^{-1}hg_i$. It follows that
$$
\deg w=g_i^{-1}hg_i g_i^{-1}g_j=g_i^{-1}hg_j=\deg(E_{\mu\nu}\otimes x_h)=
\deg\varphi(w)
$$
i.e. $\varphi$ is an isomorphism of graded algebras, and the proof of
Theorem \ref{T3} is complete.
\hfill $\Box$

{\em Remark.} The restriction on the characteristic of $F$ is used in the 
proof  of Theorem \ref{T3} only to make sure that any graded simple 
algebra $R$ is semiprimitive. Actually Theorem \ref{T3} remains true if we discard the 
restriction on the characteristic of $F$ but require instead that $R$ is semiprimitive.

\begin{Corollary}\label{C1}
Let $F$ be an algebraically closed field and $G$ a finite group such 
that $|G|^{-1}\in F$. Then
\begin{enumerate}
\item[\rm (i)] The number of pairwise non-isomorphic finite-dimensional $G$-graded division 
algebras is finite;

\item[\rm (ii)] Given an integer $n\ge 1$, the number of $n$-dimensional $G$-graded 
simple algebras is finite.
\end{enumerate}
\end{Corollary}

{\em Proof.} \begin{enumerate}
\item[\rm (i)] By Theorem \ref{T2} any graded division algebra is 
isomorphic to a twisted group ring $F^\sigma[H]$ where $H$ is a subgroup of 
$G$ and $\sigma$ is a 2-cocycle on $H$. It is known (see, for example, 
\cite[Section 1]{P1}) that two cocycles $\sigma$ and $\tau$ define the same 
twisted group algebra if there exists a map $\delta: H\rightarrow F^\ast$ 
such that
$$
\tau(x,y)=\delta(xy)\delta^{-1}(x)\delta^{-1}(y)\sigma(x,y).
$$
In particular,  any $2$-coboundary $\gamma$, i.e. a $2$-cocycle of the type $\gamma(x,y)=
\delta(xy)\delta^{-1}(x)\delta^{-1}(y)$ always defines an ordinary group ring. All cocycles form an abelian group $A$, all coboundaries form
a subgroup $B$ of $A$, and the factor-group $A/B$ is the second 
cohomology group $H^2(H,F^\ast)$, which is known to be finite for $H$ finite \cite{}. Therefore the number of non-isomorphic algebras of the type
$F^\sigma[H]$ is finite. Since the number of subgroups $H\subset G$ is finite, we get the first statement.

\item[\rm (ii)] By Theorem \ref{T3} any $G$-graded simple finite-dimensional algebra is
isomorphic to $M_k(C)$ where $C$ is a graded division algebra and the grading on $M_k(F)$ is defined by a $k$-tuple $(g_1,\ldots,g_k)\in G^k$.
Since the
number of $k$-tuples is finite, we can use part (i) of the corollary to complete the proof.\epf
\end{enumerate}

In the paper \cite{E1} the author proved that given an abelian group $G$, a 
$G$-graded ring $R$ is simple if and only if $R$ is graded simple and the
center of $R$ is a field. In \cite{E2} he extends this result to an arbitrary
hypercentral group. Clearly, $M_k(F^\sigma[H])$ is simple if and only if the 
twisted group algebra $F^\sigma[H]$ is simple and any finite-dimensional 
semiprimitive algebra over $F$ is simple if and only if its center equals 
$F$. The following result is thus an extension of the results of \cite{E1},\cite{E2}.

\begin{Corollary}\label{C2}
Let $R=\oplus_{g\in G}R_g$ be a finite-dimensional $G$-graded algebra
over an algebraically closed field $F$ such that either ${\rm char}~F=0$ or 
${\rm char}~F$ is coprime to the order of any finite subgroup of $G$.
Then $R$ is simple if and only if $R$ is graded simple and the
center of $R$ coincides with $F$.\epf
\end{Corollary}

It is known that a twisted group algebra of a finite group is semiprimitive 
under our restriction on the characteristic of $F$. In the case of abelian groups we 
can also say that all simple components are isomorphic.

\begin{Corollary}\label{C3}
Let $G$ be a finite abelian group and $F$ an algebraically closed 
field such that $|G|^{-1}\in F$. Then any twisted group algebra $F^\sigma[G]$
is the direct sum of pairwise isomorphic simple ideals.
\end{Corollary}


\begin{thebibliography}{99}

\bibitem{BaShest}  Bahturin, Yuri; Shestakov, Ivan, {\it Gradings of simple 
 Jordan algebras and their relation to the gradings of simple associative 
 algebras}, Special issue dedicated to Alexei Ivanovich Kostrikin. Comm. 
 Algebra 29 (2001), no. 9, 4095--4102.
 
 \bibitem{BSh} Bahturin, Y., Shestakov, I., {\it Group gradings on Jordan algebras $M_n^{(+)}$}, to be submitted.
 
 \bibitem{BShZ}Y. Bahturin, I. Shestakov, M. Zaicev, {\it Gradings
  on Simple Jordan and Lie  Algebras}, J. Algebra, {\bf 283} (2005), no.~2, 849--868.
  
  \bibitem{BSegZ}  Bahturin, Yu. A.; Sehgal, S. K.; Zaicev, M. V. {\it Group 
 gradings on associative algebras}, J. Algebra 241 (2001), no. 2, 677--698.
 
 \bibitem{BZ}  Bahturin, Yu. A.; Zaicev, M. V. {\it Group gradings on matrix 
 algebras},  Canad. Math. Bull. 45 (2002), no. 4, 499--508.
 
 \bibitem{BZA} Bahturin, Y., Zaicev, M., {\it Group gradings on simple Lie algebras of type ``A''}, Preprint Dept. Math., University of Palermo, 275(2005), 30 pp.
 
 
 \bibitem{CM}  Cohen, M.; Montgomery, S., {\it Group-graded rings, smash products, 
 and group actions}, Trans. Amer. Math. Soc. 282 (1984), no. 1, 237--258. 
 
  
 \bibitem{DINR}  D\u asc\u alescu, S.; Ion, B.; N\u ast\u asescu, C.; Rios 
 Montes, J. {\it Group gradings on full matrix rings},  
 J. Algebra 220 (1999), no. 2, 709--728.
 
 
 
 \bibitem{E}  Elduque, Alberto,  {\it Gradings on octonions}, J. Algebra 207 (1998), 
 no. 1, 342--354.
 
       
 \bibitem{H} Herstein, I. N. \textsc{Noncommutative rings}, Carus Mathematical Monographs, No. 15, John Wiley \& Sons, Inc., New York, 1968, xi+199 pp.
  
  \bibitem{J}  Jacobson, Nathan, \textsc{The Theory of Rings}, American Mathematical 
 Society Mathematical Surveys, vol. I. American Mathematical Society, New 
 York, 1943. vi+150 pp. 
 
 \bibitem{E1}  Jespers, Eric, {\it Simple abelian-group graded rings}, Boll. Un. 
               Mat. Ital. A (7) 3 (1989), no. 1, 103--106.
 
 \bibitem{E2}  Jespers, Eric, {\it Simple graded rings}, Comm. Algebra 21 (1993), no. 7, 2437--2444.   
  
  \bibitem{Kac}
  Kac, V. G. {\it Graded Lie algebras and symmetric spaces}, (Russian), Funkcional. 
  Anal. i Prilo\v zen. 2 1968 no. 2, 93--94.
 
 \bibitem{K} Kantor, I.L. {\it Certain generalizations of Jordan algebras}, (Russian), Trudy 
 Sem. Vektor. Tenzor. Anal. 16 (1972), 407--499.
 
 
 \bibitem{NVO}  N\u ast\u asescu, C.; van Oystaeyen, F., \textsc{Graded ring theory}, 
 North-Holland Mathematical Library, 28. North-Holland Publishing Co., 
 Amsterdam-New York, 1982. ix+340 pp.
 
 \bibitem{P1} Passman, Donald S., \textsc{The algebraic structure of group rings},
 Pure and Applied Mathematics. Wiley-Interscience [John Wiley \& Sons], New 
 York-London-Sydney, 1977. xiv+720 pp. 
 
 \bibitem{P2}  Passman, Donald S. \textsc{Infinite crossed products}, Pure and 
 Applied Mathematics, 135. Academic Press, Inc., Boston, MA, 1989. xii+468 
 pp.        
 
 
 \bibitem{Smirn}  Smirnov, Oleg N., {\it Simple associative algebras with finite 
  $\bf Z$-grading}, J. Algebra 196 (1997), no. 1, 171--184.
 
 
 \bibitem{ZS} Za\u\i tsev, M. V.; Segal, S. K., {\it Finite gradings of simple 
 Artinian rings}, (Russian), Vestnik Moskov. Univ. Ser. I Mat. Mekh. 2001, , 
 no. 3, 21--24, 77; translation in Moscow Univ. Math. Bull. 56 (2001), no. 
 3, 21--24
 
 
 \bibitem{W} Wall, C. T. C. 
 {\it Graded Brauer groups}, J. Reine Angew. Math. 213 1963/1964 187--199.
 
 

 \bibitem{S}  Sehgal, Sudarshan K., 
 \textsc{Topics in group rings}, 
 Monographs and Textbooks in Pure and Applied Math., 50. 
 Marcel Dekker, Inc., New York, 1978. vi+251 pp. 
 
 
 
 
 \end{thebibliography}
 \end{document}